\documentclass[12pt,francais,nothms]{aclart}
\usepackage{mathptm}
\usepackage[all]{xy}

\newtheorem{thm}{Th\'eor\`eme}
\newtheorem{propp}{Proposition}

\usepackage[all]{xy}

\RequirePackage[T1]{fontenc}
\RequirePackage{mathrsfs,amsmath,amssymb}
\let\mathcal\mathscr
 
\RequirePackage{url}

\RequirePackage{xspace}
\RequirePackage[english,frenchb]{babel}

\theoremstyle{remark}

\def\ord{{\rm ord}}
\def\ac{{\overline{\rm ac}}}

\def\11{{\mathbf 1}}

\def\Def{{\rm Def}}

\def\RDef{{\rm RDef}}

\def\LPas{\cL_{\rm DP}}

\def\AA{{\mathbf A}}

\def\LL{{\mathbf L}}

\def\NN{{\mathbf N}}

\def\RR{{\mathbf R}}

\def\ZZ{{\mathbf Z}}

\def\cC{{\mathcal C}}

\def\cL{{\mathcal L}}

\def\cP{{\mathcal P}}

\mathchardef\alphag="7C0B
\mathchardef\betag="7C0C
\mathchardef\gammag="7C0D
\mathchardef\deltag="7C0E
\mathchardef\varepsilong="7C22
\mathchardef\varphig="7C27
\mathchardef\psig="7C20
\mathchardef\zetag="7C10
\mathchardef\epsilong="7C0F
\mathchardef\rhog="7C1A
\mathchardef\taug="7C1C
\mathchardef\upsilong="7C1D
\mathchardef\iotag="7C13
\mathchardef\thetag="7C12
\mathchardef\pig="7C19
\mathchardef\sigmag="7C1B
\mathchardef\etag="7C11
\mathchardef\omegag="7C21
\mathchardef\kappag="7C14
\mathchardef\lambdag="7C15
\mathchardef\mug="7C16
\mathchardef\xig="7C18
\mathchardef\chig="7C1F
\mathchardef\nug="7C17
\mathchardef\varthetag="7C23
\mathchardef\varpig="7C24
\mathchardef\varrhog="7C25
\mathchardef\varsigmag="7C26
\mathchardef\Omegag="7C0A
\mathchardef\Thetag="7C02
\mathchardef\Sigmag="7C06
\mathchardef\Deltag="7C01
\mathchardef\Phig="7C08
\mathchardef\Gammag="7C00
\mathchardef\Psig="7C09
\mathchardef\Lambdag="7C03
\mathchardef\Xig="7C04
\mathchardef\Pig="7C05
\mathchardef\Upsilong="7C07

\begin{document}

\selectlanguage{french}
\title[Fonctions constructibles et int\'egration motivique I]
  {{\normalsize\bfseries G\'eom\'etrie alg\'ebrique/\itshape Algebraic geometry \\}
  
Fonctions constructibles et int\'egration motivique I}
\alttitle{Constructible functions and motivic integration I}

\author{Raf Cluckers}
\address{Katholieke Universiteit Leuven, Department of Mathematics,
Celestijnenlaan 200B, 3001 Leu\-ven, Bel\-gium }
\email{raf.cluckers@wis.kuleuven.ac.be}
\author{Fran\c cois Loeser}

\address{{\'E}cole Normale Sup{\'e}rieure,
D{\'e}partement de math{\'e}matiques et applications,
45 rue d'Ulm,
75230 Paris Cedex 05, France
(UMR 8553 du CNRS)}
\email{Francois.Loeser@ens.fr}
\date{\today}

\begin{altabstract}We introduce a direct image formalism
for constructible motivic functions. One deduces a very
general version of motivic integration for which 
a change of variables theorem is proved.
These constructions are generalized
to the relative framework,  in which we
develop a relative version of motivic integration.
Details of constructions and proofs will
be given in \cite{cl}.
\end{altabstract}

\begin{abstract}
On introduit un formalisme d'images directes 
pour les fonctions constructibles motiviques.
On en tire une version tr\`es g\'en\'erale de l'int\'egration motivique
pour laquelle un th\'eor\`eme de changement de variables est \'etabli.
Ces constructions admettent une g\'en\'eralisation au cadre relatif, ce
qui 
permet \'egalement de d\'evelopper
une version relative de 
l'int\'egration motivique.
Les d\'etails des constructions et des preuves seront donn\'es
dans \cite{cl}.

\end{abstract}

\maketitle

\selectlanguage{french}

\section{Pr\'eliminaires}

\subsection{Langage de Denef-Pas}Dans ce travail on
fixe un corps $k$ de caract\'eristique
z\'ero et on consid\`ere des corps $K$ contenant $k$ ainsi que le 
corps des s\'eries de Laurent $K ((t))$
muni
de la valuation naturelle  $\ord : K ((t))^{\times} \rightarrow \ZZ$.
Pour $x$ dans $K ((t))$ on pose $\ac (x) = x t^{-\ord (x)} \mod t$ si $x \not=0$
et $\ac (0) = 0$.
On utilise le langage de Denef-Pas
$\LPas$.
Il s'agit d'un langage \`a trois sortes
$(\LL_{\rm Val},\LL_{\rm Res},\LL_{\rm Ord},\ord,\ac)$.
Pour la sorte
de type ${\rm Val}$, on prendra comme langage
$\LL_{\rm Val}$ le langage des anneaux $\LL_{\rm
Rings}=(+,-,\cdot,0,1)$,
de m\^eme, on prendra le langage
$\LL_{\rm Res}=\LL_{\rm Rings}$ pour la sorte de type ${\rm
Res}$,
tandis que pour la sorte de type 
${\rm Ord}$
on prendra le langage de Presburger
$$
\LL_{\rm PR} = \{+, 0, 1, \leq\} \cup \{\equiv_n\ \mid n\in \NN,\
n
>1\},
$$
avec  $\equiv_n$ la relation d'\'equivalence modulo $n$.
Les symboles $\ord$ et $\ac$ 
seront interpr\'et\'es respectivement comme la valuation 
et la composante angulaire.
Ainsi $(K ((t)),K,\ZZ)$ est une structure pour $\LPas$.
Les formules du premier ordre
dans le langage $\LPas$ sont construites \`a partir des symboles
de $\LPas$ ainsi que de variables,
des connecteurs logiques $\wedge$, $\vee$, $\neg$,
des quanteurs  $\exists$, $\forall$ 
et du symbole de l'\'egalit\'e $=$.
En g\'en\'eral on dispose \'egalement d'un symbole de
constante dans la sorte
de type ${\rm Val}$, resp. {\rm Res}, pour tout \'el\'ement de $k ((t))$
resp. $k$.

\subsection{Sous-assignements}Soit $F : \cC \rightarrow {\rm Ens}$ un foncteur
d'une cat\'egorie $\cC$ \`a valeurs dans celle des ensembles.
Rappelons qu'un sous-assignement $h$ de $F$ est la donn\'ee,
pour tout objet $C$ de $\cC$, d'un sous-ensemble
$h (C)$ de $F (C)$, cf. \cite{JAMS}. 
Les op\'erations et notations usuelles de la th\'eorie des ensembles
s'\'etendent trivialement aux sous-assignements. Ainsi pour deux 
sous-assignements $h$ et $h'$ d'un m\^eme foncteur, on d\'efinit
des sous-assignements
$h \cup h'$, $h \cap h'$ et la relation $h \subset h'$, etc.
Si $h \subset h'$ on dit que $h$ est un sous-assignement de $h'$.
Un morphisme $f : h \rightarrow h'$
entre sous-assignements de foncteurs $F_1$ et $F_2$ est la donn\'ee
pour tout $C$ d'une application $f (C) : h (C) \rightarrow h'(C)$.
On d\'efinit ais\'ement
le sous-assignement $f (h)$ de $F_2$ ainsi que le graphe
de $f$, un sous-assignement de
$F_1 \times F_2$, cf. \cite{JAMS}.

\subsection{Sous-assignements
d\'efinissables}Soit $F_k$ la cat\'egorie des corps contenant $k$.
On note $h [m, n, r]$ le foncteur $F_k \rightarrow {\rm Ens}$
donn\'e par $h [m, n, r] (K) = K ((t))^m \times K^n \times \ZZ^r$.
A toute formule $\varphi$ dans $\LPas$
\`a coefficients dans $k ((t))$, resp. $k$, dans la sorte 
valu\'ee, resp. r\'esiduelle, ayant respectivement $m$, $n$ et $r$ 
variables libres dans les diff\'erents types, on associe 
un sous-assignement $h_{\varphi}$ de $h [m, n, r]$,
en prenant pour  $h_{\varphi} (K)$ le sous-ensemble de
$h [m, n, r] (K)$ form\'e des points satisfaisant $\varphi$.
Un tel sous-assignement sera appel\'e d\'efinissable.
On d\'efinit une cat\'egorie $\Def_k$
en prenant comme objets les 
sous-assignements d\'efinissables d'un
$h [m, n, r]$. Les morphismes dans $\Def_k$
sont les morphismes $f : h \rightarrow h'$
avec $h$ et $h'$ sous-assignements d\'efinissables de
$h [m, n, r]$ et $h [m', n', r']$ respectivement
dont le graphe est d\'efinissable.
Si $S$ est un objet de $\Def_k$, on note $\Def_S$ la cat\'egorie
des morphismes $X \rightarrow S$ dans $\Def_k$.
On \'ecrit $S [m, n, r]$ pour $S \times h  [m, n, r]$.
Un point $x$ de $S$ est un couple $(x_0, K)$ avec
$K$ dans $F_k$ et $x_0$ un point de $S (K)$. 
On note $\vert S \vert$ l'ensemble des points de $S$.
On pose alors
$k (x) = K$ et on dispose d'un foncteur ``fibre en $x$''
$i_x^* : \Def_S \rightarrow \Def_{k (x)}$.

\subsection{Dimension}A toute  sous-vari\'et\'e alg\'ebrique $Z$ de
$\AA^m_{k ((t))}$ on associe le sous-assignement d\'efinissable
$h_Z$ de $h [m, 0, 0]$ donn\'e par $h_Z (K) = Z (K ((t)))$.
L'adh\'erence de Zariski d'un sous-assignement $S$ de 
$h [m, 0, 0]$ est l'intersection $W$ des 
sous-vari\'et\'es alg\'ebriques $Z$ de 
$\AA^m_{k ((t))}$ telles que $S \subset h_Z$. On d\'efinit la dimension
de $S$ comme
$\dim S := \dim W$. Plus g\'en\'eralement, si $S$ est un
sous-assignement
de $h_Z$ of $h [m, n, r]$,
on d\'efinit $\dim S$ comme la dimension
de l'image de $S$ par la projection 
$h [m, n, r] \rightarrow h [m, 0, 0]$.
La d\'emonstration de l'\'enonc\'e suivant n'est pas triviale 
et repose sur des r\'esultats de Denef et Pas \cite{Pas}
et van den Dries \cite{vdDries}.

\begin{propp}Deux objets isomorphes de $\Def_k$
ont m\^eme dimension.
\end{propp}

\section{Fonctions constructibles}

\subsection{}Fixons $S$ dans
$\Def_k$.
On consid\`ere la sous-cat\'egorie $\RDef_S$ de
$\Def_S$ form\'ee des sous-assignements d\'efinissables $Z$
de $S \times h [0, n, 0]$, pour $n$ variable, le morphisme
$Z \rightarrow S$ \'etant induit par la projection sur $S$.
On note $SK_0 (\RDef_S)$ le quotient du
semi-groupe libre sur les classes
d'isomorphisme d'objets $[Z \rightarrow S]$
de 
$\RDef_S$  par les relations
$[\emptyset \rightarrow S] = 0$ et
$[(Y \cup Y') \rightarrow S] + [(Y \cap Y') \rightarrow S]
= [Y \rightarrow S] + [Y' \rightarrow
S]$ et $K_0 (\RDef_S)$ le groupe ab\'elien associ\'e.
Noter que le morphisme
$SK_0 (\RDef_S) \rightarrow K_0 (\RDef_S)$
n'est pas injectif.
Le produit cart\'esien induit une unique
structure de semi-anneau sur
$SK_0 (\RDef_S)$ et d'anneau
sur $K_0 (\RDef_S)$. 
Pour tout  morphisme $f : S \rightarrow S'$
dans $\Def_k$, on dispose d'un morphisme
$f^* : SK_0 (\RDef_{S'}) \rightarrow
SK_0 (\RDef_S)$ induit par produit fibr\'e.
Si $f : S \rightarrow S'$
est un objet dans $\RDef_{S'}$,
la composition avec $f$ induit 
un morphisme $f_! : SK_0 (\RDef_S) \rightarrow
SK_0 (\RDef_{S'})$. Ces constructions s'\'etendent \`a $K_0$.
On consid\`ere
l'anneau $A = \ZZ [\LL, \LL^{-1}, \Bigl(\frac{1}{1 - \LL^{-i}}\Bigr)_{i > 0}]$.
Pour tout r\'eel $q > 1$ on note $\vartheta_q : A \rightarrow \RR$
le morphisme envoyant $\LL$ sur $q$. On note $A_+$
le sous-semi-groupe de $A$ form\'e
des $a$ tels que $\vartheta_q (a) \geq 0$
pour tout $q > 1$.
On note $\cP (S)$ le sous-anneau de l'anneau des fonctions
$\vert S \vert \rightarrow A$ engendr\'e par les constantes,
les fonctions d\'efinissables $S \rightarrow \ZZ$ et les
fonctions de la forme $\LL^{\beta}$ avec $\beta$
d\'efinissable $S \rightarrow \ZZ$. On note $\cP_+ (S)$ le semi-anneau
form\'e des fonctions de $\cP (S)$ \`a valeurs dans $A_+$.

\subsection{}\label{fc}
Si $Y$ est un sous-assignement d\'efinissable de $S$, on note
${\bf 1}_Y$ la fonction de $\cP (S)$ valant $1$ sur $Y$ et $0$
ailleurs. On note $\cP^0 (S)$, resp.
$\cP^0_+ (S)$,
le sous-anneau de $\cP (S) $,
resp. le sous-semi-anneau de $\cP_+ (S) $,
engendr\'e par de telles fonctions et par la fonction constante
$\LL - 1$.
On note $\LL$ et $\LL - 1$ la classe
de $S[0,1,0]$ et $S \times h_{\AA^1_k \setminus \{0\}}$ dans
$SK_0 (\RDef_S)$ et dans $K_0 (\RDef_S)$. 
On a des morphismes naturels 
$\cP^0 (S) \rightarrow K_0 (\RDef_S)$
et
$\cP^0_+ (S) \rightarrow SK_0 (\RDef_S)$
envoyant ${\bf 1}_Y$ sur $[Y \rightarrow S]$ et $\LL - 1$ sur $\LL - 1$.
Finalement on d\'efinit le semi-anneau des fonctions
constructibles positives par
$\cC_+ (S) = SK_0 (\RDef_S) \otimes_{\cP^0_+ (S)} \cP_+ (S)$
et l'anneau des
fonctions
constructibles par
$\cC (S) = K_0 (\RDef_S) \otimes_{\cP^0 (S)} \cP (S)$.
Si $f : S \rightarrow S'$ est un morphisme dans
$\Def_k$ le morphisme $f^*$ a une extension naturelle
en $f^* : \cC_+ (S') \rightarrow \cC_+ (S)$.
Si, de plus, $f$ est un morphisme
dans $\RDef_{S'}$, le morphisme $f_!$ admet
une extension naturelle en 
$f_! : \cC_+ (S) \rightarrow \cC_+(S')$.
Ces constructions s'\'etendent \`a $\cC$.

\subsection{}\label{dd}
Soit $\varphi$ une fonction dans
$\cP (S [0, 0, r])$.
On dit que $\varphi$ est 
$S$-int\'egrable si pour tout
$q > 1$ et tout $x$ dans $\vert S\vert$
la s\'erie $\sum_{i \in \ZZ^r} \vartheta_q (\varphi (x, i))$
est sommable. On d\'emontre 
que si $\varphi$ est 
$S$-int\'egrable il existe une unique fonction $\mu_S (\varphi)$
dans $\cP (S)$ telle que
$\vartheta_q (\mu_S (\varphi) (x))$ soit \'egale \`a la somme de la s\'erie
pr\'ec\'edente pour tout
$q > 1$ et tout $x$ dans $\vert S\vert$.
On note ${\rm I}_S \cP_+ (S[0, 0, r])$
l'ensemble des fonctions $S$-int\'egrables
dans $ \cP_+ (S[0, 0, r])$
et on pose
${\rm I}_S \cC_+ (S[0, 0, r]) = \cC_+ (S)  \otimes_{\cP_+ (S)} {\rm I}_S \cP_+ (S[0, 0, r])$.
C'est un sous $\cC_+ (S)$-semi-module de
$ \cC_+ (S[0, 0, r])$ et $\mu_S$ s'\'etend par tensorisation en un morphisme
$\mu_S : {\rm I}_S \cC_+ (S[0, 0, r]) \rightarrow \cC_+ (S)$.

\subsection{}Pour tout entier
$d$, on note $\cC_+^{\leq d} (S)$ l'id\'eal de $\cC_+ (S)$ engendr\'e
par les fonctions $\11_{Z}$ avec $Z$ sous-assignement d\'efinissable de $S$
et $\dim Z \leq d$.
On pose $C_+ (S) = \oplus_d  C^d_+ (S)$ avec 
$C^d_+ (S) := \cC_+^{\leq d} (S) / \cC_+^{\leq d-1} (S)$.
C'est un semi-groupe ab\'elien gradu\'e, ainsi qu'un
$\cC (S)_+$-semi-module. Ses \'el\'ements sont les Fonctions
constructibles positives sur $S$. Si $\varphi$
est une fonction appartenant \`a $\cC_+^{\leq d} (S)$
mais pas \`a $\cC_+^{\leq d - 1} (S)$
on note $[\varphi]$ son image dans $C^d_+ (S)$.
Si $c : S \rightarrow h [1, 0, 0]$ est un morphisme,
resp. $f : S \rightarrow S'$ est un isomorphisme, dans $\Def_k$,
on d\'emontre que la  fonction ordre du gradient $\ord {\rm grad} \, c$,
resp. la fonction ordre du jacobien $\ord {\rm jac} f$,
qui n'est d\'efinie que presque partout, est \'egale presque partout \`a une fonction
d\'efinissable, et en particulier on peut d\'efinir 
$\LL^{\sup (0, - \ord {\rm grad} \, c)}$
et 
$\LL^{- \ord {\rm jac} f}$
dans $C_+^d (S)$, pour $S$ de dimension $d$. On d\'efinit de m\^eme $C (S)$
\`a partir de $\cC (S)$.

\section{Int\'egration motivique : le r\'esultat principal}
\begin{thm}\label{mt}Soit $k$ un corps de caract\'eristique z\'ero
et soit $S$ dans $\Def_k$.
Il existe un unique foncteur $Z \mapsto
{\rm I}_S C_+(Z)$ de $\Def_S$ dans la cat\'egorie
des semi-groupes ab\'eliens, le foncteur des Fonctions $S$-int\'egrables,
associant \`a tout morphisme $f : Z \rightarrow Y$ dans
$\Def_S$ un morphisme $f_! : {\rm I}_S C_+(Z) \rightarrow {\rm I}_S C_+(Y)$ 
tel que:

{\rm (A0)} Pour tout $Z$ dans $\Def_S$, ${\rm I}_S C_+(Z)$ 
est un sous-semi-groupe gradu\'e de $C_+ (Z)$ ; ${\rm I}_S C_+(S) = C_+ (S)$. 

{\rm (A1a)} Si $S \rightarrow S'$ est un morphisme dans $\Def_k$
et $Z$ est dans
$\Def_S$, alors ${\rm I}_{S'} C_+(Z) \subset {\rm I}_{S} C_+(Z)$, et pour
$\varphi$ dans ${\rm I}_{S'} C_+(Z)$, $f_! (\varphi)$ est le m\^eme,
consid\'er\'e
dans ${\rm I}_{S'}$ ou dans ${\rm I}_{S}$.

{\rm (A1b)} Une Fonction positive $\varphi$ sur $Z$ est $S$-int\'egrable si
et seulement si elle est $Y$-int\'egrable et $f_! (\varphi)$ est
$S$-int\'egrable.

{\rm (A2)} Si $Z$ est la r\'eunion disjointe de deux sous-assignements 
d\'efinissables $Z_1$ et $Z_2$, alors l'isomorphisme
$C_+ (Z) \simeq C_+ (Z_1) \oplus C_+ (Z_2)$ induit un isomorphisme
${\rm I}_S C_+ (Z) \simeq {\rm I}_S C_+ (Z_1) \oplus {\rm I}_S C_+ (Z_2)$,
sous lequel on a $f_! = f_{|Z_1 !} \oplus f_{|Z_2 !}$.

{\rm (A3)} Pour tout $\alpha$ dans $\cC_+ (Y)$ et tout $\beta$ dans 
${\rm I}_S C_+ (Z)$, 
$\alpha f_! (\beta)$ est $S$-int\'egrable si et seulement si
$f^* (\alpha) \beta$ l'est et dans ce cas
$f_! (  f^* (\alpha) \beta)  = \alpha f_! (\beta)$.

{\rm (A4)} Si $i : Z \hookrightarrow Z'$ est l'inclusion de
sous-assignements 
d\'efinissables d'un m\^eme objet de $\Def_S$, $i_!$ est induit par
le prolongement par z\'ero en dehors de $Z$ et envoie injectivement
${\rm I}_S C_+ (Z)$ dans ${\rm I}_S C_+ (Z')$.

{\rm (A5)} Soit $Y$ dans $\Def_S$ et soit $\pi$ la projection
$Y [0, n, 0] \rightarrow Y$. Une Fonction
$[\varphi]$ dans $C_+ (Y [0, n, 0])$ est $S$-int\'egrable si et seulement si
$[\pi_! (\varphi)]$ l'est (avec les notations de \ref{fc}) et dans ce cas
$\pi_! ([\varphi]) = [\pi_! (\varphi)]$.

{\rm (A6)}  Soit $Y$ dans $\Def_S$ et soit $\pi$ la projection
$Y [0, 0, r] \rightarrow Y$. Une Fonction
$[\varphi]$ dans $C_+ (Y [0, 0, r])$ est $S$-int\'egrable si et seulement si
il existe $\varphi'$ avec $[\varphi'] = [\varphi]$ qui soit $Y$-int\'egrable au sens de
\ref{dd} et telle que $[\mu_Y (\varphi')]$ soit $S$-int\'egrable. On pose alors
$\pi_! ([\varphi]) = [\mu_Y (\varphi')]$.

{\rm (A7)}  Soit $Y$ dans $\Def_S$ et soit $Z$ le sous-assignement de
$Y [1, 0, 0] $ d\'efini par $\ord (z - c (y)) = \alpha (y)$
et $\ac (z - c (y)) = \xi (y)$, avec $z$ la coordonn\'ee sur le facteur
$\AA^1_{k ((t))}$
et $\alpha, \xi, c$ des fonctions d\'efinissables
sur $Y$ respectivement \`a valeurs dans $\ZZ$, $h[0, 1, 0] \setminus \{0\}$ et $h[1, 0, 0]$.
On consid\`ere $f : Z \rightarrow Y$ induit par la projection.
Alors $[\11_Z]$ est $S$-int\'egrable si et seulement si $\LL^{-\alpha - 1} [\11_Y]$
l'est et dans ce cas $f_! ([\11_Z]) = \LL^{-\alpha - 1} [\11_Y]$.

{\rm (A8)} Soit $Y$ dans $\Def_S$ et soit $Z$ le sous-assignement de
$Y [1, 0, 0] $ d\'efini par $z - c (y) = 0$
avec $z$ la coordonn\'ee sur le facteur $\AA^1_{k ((t))}$
et $c$ un morphisme $Y \rightarrow h[1, 0, 0]$.
On consid\`ere $f : Z \rightarrow Y$ induit par la projection.
Alors $[\11_Z]$ est $S$-int\'egrable si et seulement si 
$\LL^{\sup (0, - \ord {\rm grad} c)}$
l'est et dans ce cas $f_! ([\11_Z]) = \LL^{\sup (0, - \ord {\rm grad} c)}$.
\end{thm}

La preuve du th\'eor\`eme  utilise de fa\c con essentielle le th\'eor\`eme
de d\'ecomposition en cellules de Denef et Pas \cite{Pas}.

On d\'efinit  ${\rm I}_S C (Y)$ comme le sous-groupe
de $C (Y)$ engendr\'e par l'image de ${\rm I}_S C_+ (Y)$.
On d\'emontre que si $f : Y \rightarrow Y'$
est un morphisme dans 
$\Def_S$, le morphisme
$f_! : {\rm I}_S C_+(Y) \rightarrow {\rm I}_S C_+(Y')$
admet une extension naturelle en
$f_! : {\rm I}_S C(Y) \rightarrow {\rm I}_S C(Y')$.

L'\'enonc\'e suivant est une version g\'en\'erale des
th\'eor\`emes de changement de variable
de \cite{arcs} et \cite{JAMS}.

\begin{thm}Soit $f : X \rightarrow Y$ un isomorphisme entre sous-assignements
d\'efinissables de dimension $d$. Pour toute fonction $\varphi$
dans $\cC^{\leq d} _+(Y)$ ayant une classe non nulle
dans $C^d (Y)_+$, $[f^* (\varphi)]$ est $Y$-int\'egrable
et
$f_! [f^* (\varphi)] = \LL^{\ord {\rm jac} f \circ f^{-1}} [\varphi]$.
On a un \'enonc\'e similaire dans $C$.
\end{thm}

Quand $S$ est \'egal \`a $h [0, 0, 0]$,
i.e. \`a l'objet final de $\Def_k$, on \'ecrit ${\rm I} C_+ (Z)$ pour
 ${\rm I}_S C_+ (Z)$ et on dira int\'egrable pour $S$-int\'egrable,
 de m\^eme pour $C$.
Notons que ${\rm I} C_+ (h [0, 0, 0]) = C_+ (h [0, 0, 0]) =
SK_0 (\RDef_k) \otimes_{\NN [\LL - 1]} A_+$ et 
que
${\rm I}C (h [0, 0, 0]) =  K_0 (\RDef_k) \otimes_{\ZZ [\LL]} A$.
Pour  $\phi$ dans ${\rm I} C_+ (Z)$, ou dans ${\rm I} C (Z)$, on d\'efinit l'int\'egrale motivique
$\mu (\varphi)$ par $\mu (\varphi) = f_! (\varphi)$ avec
$f$ le morphisme $Z \rightarrow h [0, 0, 0]$.
Les relations de cette nouvelle construction avec les
constructions ant\'erieures de l'int\'egration motivique, tant dans
sa version  g\'eom\'etrique, introduite dans 
\cite{K} et d\'evelopp\'ee dans \cite{arcs}, que dans sa version arithm\'etique 
\cite{JAMS}, ainsi qu'avec
l'int\'egration $p$-adique seront explicit\'ees ult\'erieurement.

\section{Int\'egrales d\'ependant d'un param\`etre}
On fixe $\Lambda$ dans $\Def_k$ qui joue le r\^ole d'un espace de param\`etres.
Pour $S$ dans $\Def_{\Lambda}$,
on consid\`ere l'id\'eal $\cC^{\leq d} (S \rightarrow \Lambda)$
de $\cC_+ (S)$ engendr\'e
par les fonctions $\11_{Z}$ avec $Z$ sous-assignement d\'efinissable de $S$
tel que toutes les fibres de $Z \rightarrow \Lambda$ soient de dimension $\leq d$.
On pose $C_+ (S \rightarrow \Lambda) = \oplus_d  C^d_+ (S \rightarrow \Lambda)$ avec 
$C^d_+ (S \rightarrow \Lambda) := \cC^{\leq d}_+ (S \rightarrow \Lambda)  / \cC_+^{\leq d-1} (S \rightarrow \Lambda)$.
C'est un semi-groupe ab\'elien gradu\'e (et aussi un
$\cC_+ (S)$-semi-module).  Si  $\varphi$
appartient \`a $\cC_+^{\leq d} (S \rightarrow \Lambda)$
mais pas \`a $\cC_+^{\leq d - 1} (S \rightarrow \Lambda)$
on note $[\varphi]$ son image dans $C_+^d (S \rightarrow \Lambda)$.
On a l'analogue suivant du th\'eor\`eme \ref{mt}.

\begin{thm}\label{mtr}Soit $k$ un corps de caract\'eristique z\'ero,
soit $\Lambda$ dans $\Def_k$ et soit $S$ dans $\Def_{\Lambda}$.
Il existe un unique foncteur $Z \mapsto
{\rm I}_S C_+(Z \rightarrow \Lambda)$ de $\Def_S$ dans la cat\'egorie
des semi-groupes ab\'eliens, 
associant \`a tout morphisme $f : Z \rightarrow Y$ dans
$\Def_S$ un morphisme $f_{! \Lambda} :
{\rm I}_S C_+(Z\rightarrow \Lambda))
\rightarrow {\rm I}_S C_+(Y\rightarrow \Lambda))$ 
v\'erifiant les analogues de (A0)-(A8) obtenus 
en rempla\c cant 
$ C_+ (\_)$ par $C_+ (\_ \rightarrow \Lambda)$ et
$\ord {\rm grad}$
par son analogue relatif
$\ord {\rm grad}_{\Lambda}$, 
pour lequel les d\'eriv\'ees partielles ne sont prises
que par rapport aux variables dans les fibres de la projection sur
$\Lambda$.
\end{thm}

Pour $f : Z \rightarrow \Lambda$ dans
$\Def_{\Lambda}$, on dispose ainsi d'un morphisme
$\mu_{\Lambda} := 
f_{!\Lambda } : {\rm I}_{\Lambda} C_+ (Z \rightarrow \Lambda) \rightarrow \cC_+ (\Lambda)
= {\rm I}_{\Lambda} C_+ (\Lambda \rightarrow \Lambda)$
qui correspond \`a
l'int\'egration dans les fibres de $\Lambda$ d'apr\`es
l'\'enonc\'e suivant.

\begin{propp}\label{ll}Soit $\varphi$
une Fonction dans $C_+(Z \rightarrow \Lambda)$.
Elle appartient \`a
${\rm I}_{\Lambda} C_+(Z \rightarrow \Lambda)$ si et seulement si
pour tout point $\lambda$ de $\Lambda$, la restriction $\varphi_{\lambda}$ 
de
$\varphi$ \`a la fibre de $Z$ en $\lambda$
est int\'egrable. L'int\'egrale motivique de
$\varphi_{\lambda}$ 
est alors \'egale
\`a $i_{\lambda}^* (\mu_{\Lambda} (\varphi))$, pour tout $\lambda$.
\end{propp}

Bien entendu on peut \'egalement d\'efinir l'analogue relatif 
$C (S \rightarrow \Lambda)$ de $C (S)$,
et \'etendre la notion d'int\'egrabilit\'e et la construction de
$f_{!\Lambda}$ \`a ce cadre.

\subsection*{}{\small Pendant la r\'ealisation de ce projet, 
le premier auteur  \'etait chercheur postdoctoral du Fonds 
de Recherche Scientifique - Flandres (Belgique)
et il a b\'en\'efici\'e du soutien partiel du projet europ\'een
EAGER.}


 \bibliographystyle{smfplain}
 \bibliography{aclab,acl}

\providecommand{\noopsort}[1]{}\providecommand{\url}[1]{\textit{#1}}
\providecommand{\bysame}{\leavevmode ---\ }
\providecommand{\og}{``}
\providecommand{\fg}{''}
\providecommand{\smfandname}{\&}
\providecommand{\smfedsname}{\'eds.}
\providecommand{\smfedname}{\'ed.}
\providecommand{\smfmastersthesisname}{M\'emoire}
\providecommand{\smfphdthesisname}{Th\`ese}
\begin{thebibliography}{10}



\bibitem{cl}{\scshape R. Cluckers {\normalfont \smfandname} F. Loeser},
\textit{Constructible motivic functions and motivic integration},
en pr\'eparation.

\bibitem{arcs}{\scshape J. Denef {\normalfont \smfandname} F. Loeser},
\textit{Germs of arcs on singular algebraic varieties
and motivic integration},
Invent. Math.
\textbf{135}
(1999),
201--232.

\bibitem{JAMS}
{\scshape J. Denef {\normalfont \smfandname} F. Loeser},
\textit{Definable sets, motives and $p$-adic integrals},
J. Amer. Math. Soc.,
\textbf{14} (2001), 429--469.




\bibitem{vdDries}
{\scshape L. van den Dries}, \textit{Dimension of definable sets, algebraic
boundedness and {H}enselian fields}, Ann. Pure Appl. Logic,
\textbf{45} (1989), 189--209.





\bibitem{K}
{\scshape M. Kontsevich}, Expos\'e \`a  Orsay,
7 d\'ecembre 1995.




\bibitem{Pas}
{\scshape J. Pas},
\textit{Uniform $p$-adic cell decomposition and local zeta functions},
J. Reine Angew. Math.,
\textbf{399}
(1989),
137--172.






















\end{thebibliography}
\end{document}